\documentclass[11pt]{article}
\usepackage{epsfig}

\def\be{\begin{equation}}
\def\ee{\end{equation}}
\def\bea{\begin{eqnarray}}
\def\eea{\end{eqnarray}}
\def\bes{\begin{eqnarray*}}
\def\ees{\end{eqnarray*}}

\def\nn{\nonumber}
\def\lb{\label}
\def\bs{\setminus}

\def\vs{{\varsigma}}

\def\R{{\bf R}}
\def\C{{\bf C}}
\def\Z{{\bf Z}}
\def\K{{\bf K}}
\def\N{{\bf N}}
\def\U{{\bf U}}

\def\Q{{\bf Q}}

\def\aa{{\alpha}}
\def\bb{{\beta}}
\def\ga{{\gamma}}

\def\th{{\theta}}

\def\om{{\omega}}
\def\Om{{\Omega}}
\def\ep{{\epsilon}}
\def\lm{{\lambda}}
\def\Lm{{\Lambda}}

\def\sg{{\sigma}}
\def\dm{{\diamond}}

\def\vf{{\varphi}}

\def\<{{\langle}}
\def\>{{\rangle}}

\def\K{{\cal K}}
\def\P{{\cal P}}

\def\rank{{\rm rank}}

\def\Sp{{\rm Sp}}

\def\ol{\overline}

\def\hb{\vrule height0.18cm width0.14cm $\,$}

\title{Non-hyperbolic closed geodesics on positively curved Finsler spheres}

\author{Huagui Duan \thanks{Partially supported by NNSF (No.11131004, 11471169), LPMC of MOE of China and
Nankai University. E-mail: duanhg@nankai.edu.cn}\\\\
School of Mathematical Sciences and LPMC\\
Nankai University\\
Tianjin 300071, The People's Republic of China\\}

\begin{document}
\date{}
\maketitle

\begin{abstract}
{\it In this paper, we prove that for every Finsler $n$-dimensional sphere $(S^n,F), n\ge 3$ with reversibility $\lambda$ and flag curvature $K$ satisfying $\left(\frac{\lambda}{1+\lambda}\right)^2<K\le 1$, there exist at least three distinct closed geodesics and at least two of them are elliptic if the number of prime closed geodesics is finite. When $n\ge 6$, these three distinct closed geodesics are non-hyperbolic.}
\end{abstract}



\renewcommand{\theequation}{\thesection.\arabic{equation}}
\renewcommand{\thefigure}{\thesection.\arabic{figure}}

\setcounter{figure}{0}
\setcounter{equation}{0}
\section{Introduction and main result}%Section 1

A closed curve on a Finsler manifold is a closed geodesic if it is
locally the shortest path connecting any two nearby points on this
curve. As usual, on any Finsler manifold
$(M, F)$, a closed geodesic $c:S^1=\R/\Z\to M$ is {\it prime}
if it is not a multiple covering (i.e., iteration) of any other
closed geodesics. Here the $m$-th iteration $c^m$ of $c$ is defined
by $c^m(t)=c(mt)$. The inverse curve $c^{-1}$ of $c$ is defined by
$c^{-1}(t)=c(1-t)$ for $t\in \R$.  Note that unlike Riemannian manifold,
the inverse curve $c^{-1}$ of a closed geodesic $c$
on a irreversible Finsler manifold need not be a geodesic.
We call two prime closed geodesics
$c$ and $d$ {\it distinct} if there is no $\th\in (0,1)$ such that
$c(t)=d(t+\th)$ for all $t\in\R$.
On a reversible Finsler (or Riemannian) manifold, two closed geodesics
$c$ and $d$ are called { \it geometrically distinct} if $
c(S^1)\neq d(S^1)$, i.e., their image sets in $M$ are distinct.
We shall omit the word {\it distinct} when we talk about more than one prime closed geodesic.

For a closed geodesic $c$ on $n$-dimensional manifold $(M,\,F)$, denote by $P_c$
the linearized Poincar\'{e} map of $c$. Then $P_c\in \Sp(2n-2)$ is symplectic.
For any $M\in \Sp(2k)$, we define the {\it elliptic height } $e(M)$
of $M$ to be the total algebraic multiplicity of all eigenvalues of
$M$ on the unit circle $\U=\{z\in\C|\; |z|=1\}$ in the complex plane
$\C$. Since $M$ is symplectic, $e(M)$ is even and $0\le e(M)\le 2k$.
A closed geodesic $c$ is called {\it elliptic} if $e(P_c)=2(n-1)$, i.e., all the
eigenvalues of $P_c$ locate on $\U$; {\it hyperbolic} if $e(P_c)=0$, i.e., all the
eigenvalues of $P_c$ locate away from $\U$;
{\it non-degenerate} if $1$ is not an eigenvalue of $P_c$. A Finsler manifold $(M,\,F)$
is called {\it bumpy} if all the closed geodesics on it are non-degenerate.

There is a famous conjecture in Riemannian geometry which claims the existence of infinitely many
closed geodesics on any compact Riemannian manifold. This conjecture
has been proved for many cases, but not yet for compact rank one symmetric
spaces except for $S^2$. The results of Franks \cite{Fra} in 1992 and Bangert \cite{Ban}
in 1993 imply that this conjecture is true for any Riemannian 2-sphere (cf. \cite{Hin2} and \cite{Hin3}).
But once one moves to the Finsler case, the conjecture
becomes false. It was quite surprising when Katok \cite{Kat}  in 1973 found
some irreversible Finsler metrics on spheres with only finitely
many closed geodesics and all closed geodesics are non-degenerate and elliptic (cf. \cite{Zil}).

Recently, index iteration theory of closed geodesics (cf. \cite{Bot} and \cite{Lon3}) has been applied to study the closed geodesic problem on Finsler manifolds. For example, Bangert and Long in \cite{BaL} show that there exist at least two closed geodesics on every $(S^2,F)$. After that, a great number of multiplicity and stability results have appeared (cf. \cite{DuL1}-\cite{DuL4}, \cite{DLW}, \cite{Lon4}, \cite{LoD}, \cite{Rad4}-\cite{Rad5}, \cite{Wan1}-\cite{Wan4} and therein).

In \cite{Rad3}, Rademacher has introduced the reversibility $\lambda=\lambda(M,F)$ of a compact Finsler manifold defined by
\bea \lambda=\max\{F(-X)\ |\ X\in TM,\ F(X)=1\}\ge 1.\nn\eea
Then Rademacher in \cite{Rad4} has obtained some results about multiplicity and the length of closed geodesics and about their stability properties. For example, let $F$ be a Finsler metric on $S^{n}$ with reversibility $\lm$ and flag curvature $K$ satisfying $\left(\frac{\lm}{1+\lm}\right)^2<K\le 1$, then there exist at least $n/2-1$ closed geodesics with length $<2n\pi$. If $\frac{9\lm^2}{4(1+\lm)^2}<K\le 1$ and $\lm<2$, then there exists a closed geodesic of elliptic-parabolic, i.e., its linearized Poincar\'{e} map split into $2$-dimensional rotations and a part whose eigenvalues are $\pm 1$. Some similar results in the Riemannian case are obtained in \cite{BTZ1} and \cite{BTZ2}.

Recently, Wang in \cite{Wan1} proved that for every Finsler $n$-dimensional sphere $S^n$ with reversibility $\lm$ and flag curvature $K$ satisfying $\left(\frac{\lm}{1+\lm}\right)^2<K\le 1$, either there exist infinitely many prime closed geodesics or there exists one elliptic closed geodesics whose linearized Poincar\'{e} map has at least one eigenvalue which is of the form $\exp(\pi i\mu)$ with an irrational $\mu$. Wang in \cite{Wan4} proved that for every Finsler $n$-dimensional sphere $S^n$ for $n\ge 6$ with reversibility $\lm$ and flag curvature $K$ satisfying $\left(\frac{\lm}{1+\lm}\right)^2<K\le 1$, either there exist infinitely many prime closed geodesics or there exists $[\frac{n}{2}]-2$ closed geodesics possessing irrational mean indices. Furthermore, assume that this metric $F$ is bumpy, in \cite{Wan2}, Wang showed that there exist $2[\frac{n+1}{2}]$ closed geodesics on $(S^n,F)$. Also in \cite{Wan2}, Wang showed that for every bumpy Finsler metric $F$ on $S^n$ satisfying $\frac{9\lm^2}{4(1+\lm)^2}<K\le 1$, there exist two prime elliptic closed geodesics provided the number of closed geodesics on $(S^n,F)$ is finite.

Very recently, the author in \cite{Dua} proved that for every Finsler $n$-dimensional sphere $(S^{n},F)$ for $n\ge 2$ with reversibility $\lm$ and flag curvature $K$ satisfying $\left(\frac{\lm}{1+\lm}\right)^2<K\le 1$, either there exist infinitely many closed geodesics, or there exist at least two elliptic closed geodesics and each linearized Poincar\'{e} map has at least one eigenvalue of the form $e^{\sqrt{-1}\th}$ with $\th$ being an irrational multiple of $\pi$.

In this paper, we generalize some above results to the following theorem.

\medskip

{\bf Theorem 1.1.} {\it For every Finsler metric $F$ on the $n$-dimensional sphere $S^n$, $n\ge 3$ with reversibility $\lm$ and flag curvature $K$ satisfying $\left(\frac{\lm}{1+\lm}\right)^2<K\le 1$, either there exist infinitely many closed geodesics, or there exist always three prime closed geodesics and at least two of them are elliptic. When $n\ge 6$, these three distinct closed geodesics are non-hyperbolic.}

\medskip

Also note that Wang in \cite{Wan1} obtained the existence of three prime closed geodesics on $(S^3,F)$ with reversibility $\lm$ and flag curvature $K$ satisfying $\left(\frac{\lm}{1+\lm}\right)^2<K\le 1$. In Section 3, we will reprove this case of $n=3$ in a more simple argument by our method. In addition, when $n\ge 10$, Theorem 1.1 is included in Theorem 1.2 of \cite{Wan4}.

Our proof of Theorem 1.1 in Section 3 contains mainly three ingredients: the common index jump theorem of \cite{LoZ}, Morse theory and some new symmetric information about index jump. In addition, we also follow some ideas from our recent preprints \cite{Dua} and \cite{DuL4}.

In this paper, let $\N$, $\N_0$, $\Z$, $\Q$, $\R$, and $\C$ denote
the sets of natural integers, non-negative integers, integers,
rational numbers, real numbers, and complex numbers respectively.
We use only singular homology modules with $\Q$-coefficients.
For an $S^1$-space $X$, we denote by $\overline{X}$ the quotient space $X/S^1$.
We define the functions
\be \left\{\matrix{[a]=\max\{k\in\Z\,|\,k\le a\}, &
           E(a)=\min\{k\in\Z\,|\,k\ge a\} , \cr
    \varphi(a)=E(a)-[a], &\{a\}=a-[a].  \cr}\right. \lb{1.1}\ee
Especially, $\varphi(a)=0$ if $ a\in\Z\,$, and $\varphi(a)=1$ if $
a\notin\Z\,$.

\setcounter{equation}{0}
\section{Morse theory and Morse indices of closed geodesics}%{Section 2}

\subsection{Morse theory for closed geodesics} % Section 2.1

Let $M=(M,F)$ be a compact Finsler manifold $(M,F)$, the space
$\Lambda=\Lambda M$ of $H^1$-maps $\gamma:S^1\rightarrow M$ has a
natural structure of Riemannian Hilbert manifolds on which the
group $S^1=\R/\Z$ acts continuously by isometries. This action is defined by
$(s\cdot\gamma)(t)=\gamma(t+s)$ for all $\gamma\in\Lm$ and $s,
t\in S^1$. For any $\gamma\in\Lambda$, the energy functional is
defined by
\be E(\gamma)=\frac{1}{2}\int_{S^1}F(\gamma(t),\dot{\gamma}(t))^2dt.
\lb{2.1}\ee
It is $C^{1,1}$ and invariant under the $S^1$-action. The
critical points of $E$ of positive energies are precisely the closed geodesics
$\gamma:S^1\to M$. The index form of the functional $E$ is well
defined along any closed geodesic $c$ on $M$, which we denote by
$E''(c)$. As usual, we denote by $i(c)$ and
$\nu(c)$ the Morse index and nullity of $E$ at $c$. In the
following, we denote by
\be \Lm^\kappa=\{d\in \Lm\;|\;E(d)\le\kappa\},\quad \Lm^{\kappa-}=\{d\in \Lm\;|\; E(d)<\kappa\},
  \quad \forall \kappa\ge 0. \nn\ee
For a closed geodesic $c$ we set $ \Lm(c)=\{\ga\in\Lm\mid E(\ga)<E(c)\}$.

Recall that respectively the mean index $\hat{i}(c)$ and the $S^1$-critical modules of $c^m$ are defined by
\be \hat{i}(c)=\lim_{m\rightarrow\infty}\frac{i(c^m)}{m}, \quad \overline{C}_*(E,c^m)
   = H_*\left((\Lm(c^m)\cup S^1\cdot c^m)/S^1,\Lm(c^m)/S^1; \Q\right).\lb{2.3}\ee

We call a closed geodesic satisfying the isolation condition, if
the following holds:

{\bf (Iso)  For all $m\in\N$ the orbit $S^1\cdot c^m$ is an
isolated critical orbit of $E$. }

Note that if the number of prime closed geodesics on a Finsler manifold
is finite, then all the closed geodesics satisfy (Iso).

If $c$ has multiplicity $m$, then the subgroup $\Z_m=\{\frac{n}{m}\mid 0\leq n<m\}$
of $S^1$ acts on $\overline{C}_*(E,c)$. As studied in p.59 of \cite{Rad2},
for all $m\in\N$, let
$H_{\ast}(X,A)^{\pm\Z_m}
   = \{[\xi]\in H_{\ast}(X,A)\,|\,T_{\ast}[\xi]=\pm [\xi]\}$,
where $T$ is a generator of the $\Z_m$-action.
On $S^1$-critical modules of $c^m$, the following lemma holds:

\medskip

{\bf Lemma 2.1.} (cf. Satz 6.11 of \cite{Rad2} and \cite{BaL})
{\it Suppose $c$ is
a prime closed geodesic on a Finsler manifold $M$ satisfying (Iso). Then
there exist $U_{c^m}^-$ and $N_{c^m}$, the so-called local negative
disk and the local characteristic manifold at $c^m$ respectively,
such that $\nu(c^m)=\dim N_{c^m}$ and
\bea \overline{C}_q( E,c^m)
&\equiv& H_q\left((\Lm(c^m)\cup S^1\cdot c^m)/S^1, \Lm(c^m)/S^1\right)\nn\\
&=& \left(H_{i(c^m)}(U_{c^m}^-\cup\{c^m\},U_{c^m}^-)
    \otimes H_{q-i(c^m)}(N_{c^m}\cup\{c^m\},N_{c^m})\right)^{+\Z_m}, \nn
\eea

(i) When $\nu(c^m)=0$, there holds
$$ \overline{C}_q( E,c^m) = \left\{\matrix{
     \Q, &\quad {\it if}\;\; i(c^m)-i(c)\in 2\Z\;\;{\it and}\;\;
                   q=i(c^m),\;  \cr
     0, &\quad {\it otherwise}, \cr}\right.  $$

(ii) When $\nu(c^m)>0$, there holds
$$ \overline{C}_q( E,c^m)=H_{q-i(c^m)}(N_{c^m}\cup\{c^m\},N_{c^m})^{\ep(c^m)\Z_m}, $$
where $\ep(c^m)=(-1)^{i(c^m)-i(c)}$.}

\medskip

Define
\be  k_j(c^m) \equiv \dim\, H_j( N_{c^m}\cup\{c^m\},N_{c^m}), \quad
     k_j^{\pm 1}(c^m) \equiv \dim\, H_j(N_{c^m}\cup\{c^m\},N_{c^m})^{\pm\Z_m}. \nn\ee

Then we have

{\bf Lemma 2.2.} (cf. \cite{Rad2}, \cite{LoD}, \cite{Wan1})
{\it Let $c$ be a  prime closed geodesic on a Finsler manifold $(M,F)$. Then

(i) For any $m\in\N$,
$k_j(c^m)=0$ for $j\not\in [0,\nu(c^m)]$.

(ii) For any $m\in\N$, $k_0(c^m)+k_{\nu(c^m)}(c^m)\le 1$ and if $k_0(c^m)+k_{\nu(c^m)}(c^m)=1$ then  $k_j(c^m)=0$ for $j\in (0,\nu(c^m))$.

(iii) For any $m\in\N$, there hold $k_0^{+1}(c^m) = k_0(c^m)$ and $k_0^{-1}(c^m) = 0$.
In particular, if $c^m$ is non-degenerate, there hold $k_0^{+1}(c^m) = k_0(c^m)=1$,
and $k_0^{-1}(c^m) = k_j^{\pm 1}(c^m)=0$ for all $j\neq 0$.

(iv) Suppose for some integer $m=np\ge 2$ with $n$ and $p\in\N$ the nullities
satisfy $\nu(c^m)=\nu(c^n)$. Then there hold $k_j(c^m)=k_j(c^n)$ and
${k}_j^{\pm 1}(c^m)={k}^{\pm 1}_j(c^n)$ for any integer $j$. }

\medskip

Let $(M,F)$ be a compact simply connected Finsler manifold with finitely many closed geodesics. It
is well known that for every prime closed geodesic $c$ on $(M,F)$, there holds either
$\hat{i}(c)>0$ and then $i(c^m)\to +\infty$ as $m\to +\infty$, or $\hat{i}(c)=0$ and
then $i(c^m)=0$ for all $m\in\N$. Denote those prime closed geodesics on $(M,F)$ with
positive mean indices by $\{c_j\}_{1\le j\le k}$.
H.-B. Rademacher in \cite{Rad1} and \cite{Rad2} established a celebrated mean index identity relating all the $c_j$s
with the global homology of $M$ (cf. Section 7, specially Satz 7.9 of \cite{Rad2}) for
compact simply connected Finsler manifolds. Here we give a brief review on this identity.

{\bf Theorem 2.3.} (Satz 7.9 of \cite{Rad2}, cf. also \cite{DuL3}, \cite{LoD} and \cite{Wan1})
{\it Assume that there exist finitely many closed geodesic on $(S^n,F)$ and
denote prime closed geodesics with positive mean indices by
$\{c_j\}_{1\le j\le k}$ for some $k\in\N$. Then the following identity holds
\be  \sum_{j=1}^k\frac{\hat{\chi}(c_j)}{\hat{i}(c_j)} = B(n,1)=\left\{\matrix{
    \frac{n+1}{2(n-1)}, &\quad n \ {\rm odd},  \cr
    -\frac{n}{2(n-1)}, &\quad n \ {\rm even}, \cr}\right.\lb{2.5}\ee
where
\be \hat{\chi}(c_j)
= \frac{1}{n(c_j)}\sum_{1\le m\le n(c_j) \atop 0\le l\le 2(n-1)}
     \chi(c_j^m)=\frac{1}{n(c_j)}\sum_{1\le m\le n(c_j) \atop 0\le l\le 2(n-1)}(-1)^{i(c_j^m)+l}k_l^{\ep(c_j^m)}(c_j^m)\in\Q     \lb{2.6}\ee
and the analytical period
$n(c_j)$ of $c_j$ is defined by (cf. \cite{LoD})
\be n(c_j) = \min\{l\in\N\,|\,\nu(c_j^l)=\max_{m\ge 1}\nu(c_j^m),\;\;
                  i(c_j^{m+l})-i(c_j^{m})\in 2\Z, \;\;\forall\,m\in\N\}. \lb{2.7}\ee}

Set $\ol{\Lm}^0=\ol{\Lambda}^0S^n =\{{\rm constant\;point\;curves\;in\;}S^n\}\cong S^n$.
Let $(X,Y)$ be a space pair such that the Betti numbers $b_i=b_i(X,Y)=\dim H_i(X,Y;\Q)$
are finite for all $i\in \Z$. As usual the {\it Poincar\'e series} of $(X,Y)$ is
defined by the formal power series $P(X, Y)=\sum_{i=0}^{\infty}b_it^i$. We need the
following well known version of results on Betti numbers and the Morse inequality.

{\bf Lemma 2.4.} (cf. Theorem 2.4 and Remark 2.5 of \cite{Rad1} and \cite{Hin1}, cf.
also Lemma 2.5 of \cite{DuL3}) {\it Let $(S^n,F)$ be a
$n$-dimensional Finsler sphere.}

(i) {\it When $n$ is odd, the Betti numbers are given by
\bea b_j
&=& \rank H_j(\Lm S^n/S^1,\Lm^0 S^n/S^1;\Q)  \nn\\
&=& \left\{\matrix{
    2,&\quad {\it if}\quad j\in \K\equiv \{k(n-1)\,|\,2\le k\in\N\},  \cr
    1,&\quad {\it if}\quad j\in \{n-1+2k\,|\,k\in\N_0\}\bs\K,  \cr
    0 &\quad {\it otherwise}. \cr}\right. \lb{2.8}\eea}

(ii) {\it When $n$ is even, the Betti numbers are given by
\bea b_j
&=& \rank H_j(\Lm S^n/S^1,\Lm^0 S^n/S^1;\Q)  \nn\\
&=& \left\{\matrix{
    2,&\quad {\it if}\quad j\in \K\equiv \{k(n-1)\,|\,3\le k\in 2\N+1\},  \cr
    1,&\quad {\it if}\quad j\in \{n-1+2k\,|\,k\in\N_0\}\bs\K,  \cr
    0 &\quad {\it otherwise}. \cr}\right.    \lb{2.9}\eea}

{\bf Theorem 2.5.} (cf. Theorem I.4.3 of \cite{Cha})
{\it Let $(M,F)$ be a Finsler manifold with finitely many closed geodesics, denoted by $\{c_j\}_{1\le j\le k}$. Set
\bea M_q =\sum_{1\le j\le k,\; m\ge 1}\dim{\ol{C}}_q(E, c^m_j),\quad q\in\Z.\nn\eea
Then for every integer $q\ge 0$ there holds }
\bea M_q - M_{q-1} + \cdots +(-1)^{q}M_0
&\ge& {b}_q - {b}_{q-1}+ \cdots + (-1)^{q}{b}_0, \lb{2.10}\\
M_q &\ge& {b}_q. \lb{2.11}\eea

\subsection{Index iteration theory of closed geodesics}%Section 2.2

In \cite{Lon1} of 1999, Y. Long established the basic normal form
decomposition of symplectic matrices. Based on this result he
further established the precise iteration formulae of indices of
symplectic paths in \cite{Lon2} of 2000. Note that this index iteration formulae works for Morse indices of iterated closed geodesics (cf. \cite{Liu} and Chap. 12 of \cite{Lon3}). Since every closed geodesic on a sphere must be orientable. Then by Theorem 1.1 of \cite{Liu}, the initial Morse index of a closed geodesic on a Finsler $S^n$ coincides with the index of a corresponding symplectic path.

As in \cite{Lon2}, denote by
\bea
N_1(\lm, b) &=& \left(\matrix{\lm & b\cr
                                0 & \lm\cr}\right), \qquad {\rm for\;}\lm=\pm 1, \; b\in\R, \lb{2.12}\\
D(\lm) &=& \left(\matrix{\lm & 0\cr
                      0 & \lm^{-1}\cr}\right), \qquad {\rm for\;}\lm\in\R\bs\{0, \pm 1\}, \lb{2.13}\\
R(\th) &=& \left(\matrix{\cos\th & -\sin\th \cr
                           \sin\th & \cos\th\cr}\right), \qquad {\rm for\;}\th\in (0,\pi)\cup (\pi,2\pi), \lb{2.14}\\
N_2(e^{\th\sqrt{-1}}, B) &=& \left(\matrix{ R(\th) & B \cr
                  0 & R(\th)\cr}\right), \qquad {\rm for\;}\th\in (0,\pi)\cup (\pi,2\pi)\;\; {\rm and}\; \nn\\
        && \qquad B=\left(\matrix{b_1 & b_2\cr
                                  b_3 & b_4\cr}\right)\; {\rm with}\; b_j\in\R, \;\;
                                         {\rm and}\;\; b_2\not= b_3. \lb{2.15}\eea
Here $N_2(e^{\th\sqrt{-1}}, B)$ is non-trivial if $(b_2-b_3)\sin\theta<0$, and trivial
if $(b_2-b_3)\sin\theta>0$.

As in \cite{Lon2}, the $\diamond$-sum (direct sum) of any two real matrices is defined by
$$ \left(\matrix{A_1 & B_1\cr C_1 & D_1\cr}\right)_{2i\times 2i}\diamond
      \left(\matrix{A_2 & B_2\cr C_2 & D_2\cr}\right)_{2j\times 2j}
=\left(\matrix{A_1 & 0 & B_1 & 0 \cr
                                   0 & A_2 & 0& B_2\cr
                                   C_1 & 0 & D_1 & 0 \cr
                                   0 & C_2 & 0 & D_2}\right). $$

For every $M\in\Sp(2n)$, the homotopy set $\Omega(M)$ of $M$ in $\Sp(2n)$ is defined by
$$ \Om(M)=\{N\in\Sp(2n)\,|\,\sg(N)\cap\U=\sg(M)\cap\U\equiv\Gamma\;\mbox{and}
                    \;\nu_{\om}(N)=\nu_{\om}(M),\ \forall\om\in\Gamma\}, $$
where $\sg(M)$ denotes the spectrum of $M$,
$\nu_{\om}(M)\equiv\dim_{\C}\ker_{\C}(M-\om I)$ for $\om\in\U$.
The component $\Om^0(M)$ of $P$ in $\Sp(2n)$ is defined by
the path connected component of $\Om(M)$ containing $M$.

\medskip

{\bf Theorem 2.6.} (cf. Theorem 7.8 of \cite{Lon1}, Theorems 1.2 and 1.3 of \cite{Lon2}, cf. also
Theorem 1.8.10, Lemma 2.3.5 and Theorem 8.3.1 of \cite{Lon3}) {\it For every $P\in\Sp(2n-2)$, there
exists a continuous path $f\in\Om^0(P)$ such that $f(0)=P$ and
\bea f(1)
&=& N_1(1,1)^{\dm p_-}\,\dm\,I_{2p_0}\,\dm\,N_1(1,-1)^{\dm p_+}
  \dm\,N_1(-1,1)^{\dm q_-}\,\dm\,(-I_{2q_0})\,\dm\,N_1(-1,-1)^{\dm q_+} \nn\\
&&\dm\,N_2(e^{\aa_{1}\sqrt{-1}},A_{1})\,\dm\,\cdots\,\dm\,N_2(e^{\aa_{r_{\ast}}\sqrt{-1}},A_{r_{\ast}})
  \dm\,N_2(e^{\bb_{1}\sqrt{-1}},B_{1})\,\dm\,\cdots\,\dm\,N_2(e^{\bb_{r_{0}}\sqrt{-1}},B_{r_{0}})\nn\\
&&\dm\,R(\th_1)\,\dm\,\cdots\,\dm\,R(\th_{r'})\,\dm\,R(\th_{r'+1})\,\dm\,\cdots\,\dm\,R(\th_r)\dm\,H(2)^{\dm h},\lb{2.16}\eea
where $\frac{\th_{j}}{2\pi}\in\Q\cap(0,1)$ for $1\le j\le r'$ and
$\frac{\th_{j}}{2\pi}\notin\Q\cap(0,1)$ for $r'+1\le j\le r$; $N_2(e^{\aa_{j}\sqrt{-1}},A_{j})$'s
are nontrivial and $N_2(e^{\bb_{j}\sqrt{-1}},B_{j})$'s are trivial, and non-negative integers
$p_-, p_0, p_+,q_-, q_0, q_+,r,r_\ast,r_0,h$ satisfy the equality
\be p_- + p_0 + p_+ + q_- + q_0 + q_+ + r + 2r_{\ast} + 2r_0 + h = n-1. \lb{2.17}\ee

Let $\ga\in\P_{\tau}(2n-2)=\{\ga\in C([0,\tau],\Sp(2n-2))\,|\,\ga(0)=I\}$. Denote the basic normal form
decomposition of $P\equiv \ga(\tau)$ by (\ref{2.16}). Then we have
\bea i(\ga^m)
&=& m(i(\ga)+p_-+p_0-r ) + 2\sum_{j=1}^r{E}\left(\frac{m\th_j}{2\pi}\right) - r   \nn\\
&&  - p_- - p_0 - {{1+(-1)^m}\over 2}(q_0+q_+)
              + 2\sum_{j=1}^{r_{\ast}}\vf\left(\frac{m\aa_j}{2\pi}\right) - 2r_{\ast}, \lb{2.18}\\
\nu(\ga^m)
 &=& \nu(\ga) + {{1+(-1)^m}\over 2}(q_-+2q_0+q_+) + 2\vs(m,\ga(\tau)),    \lb{2.19}\eea
where we denote by }
\be \vs(m,\ga(\tau)) = r - \sum_{j=1}^r\vf(\frac{m\th_j}{2\pi})
             + r_{\ast} - \sum_{j=1}^{r_{\ast}}\vf(\frac{m\aa_j}{2\pi})
             + r_0 - \sum_{j=1}^{r_0}\vf(\frac{m\bb_j}{2\pi}).    \lb{2.20}\ee

\medskip

The following is the common index jump theorem of Long and Zhu in \cite{LoZ}.

\medskip

{\bf Theorem 2.7.} (cf. Theorems 4.1-4.3 of \cite{LoZ} and \cite{Lon3}) {\it   Let $\gamma_k, k = 1,\ldots,q$ be a finite collection of symplectic paths and $M_k = \gamma_k(\tau_k)\in Sp(2n-2)$.
Suppose $\hat i(\gamma_k, 1) > 0$, for all $k = 1,\ldots,q$.
Then there exist infinitely many $(N,m_1,\ldots,m_q)\in\N^{q+1}$  such that
\bea \nu(\gamma_k, 2m_k -1)&=&  \nu(\gamma_k, 1),\lb{2.21}\\
\nu(\gamma_k, 2m_k +1)&=&   \nu(\gamma_k, 1),\lb{2.22}\\
i(\gamma_k, 2m_k -1)+\nu(\gamma_k, 2m_k -1)&=&
2N-\left(i(\gamma_k, 1)+2S^+_{M_k}(1)-\nu(\gamma_k, 1)\right),\lb{2.23}\\
i(\gamma_k, 2m_k+1)&=&2N+i(\gamma_k, 1),\lb{2.24}\\
i(\gamma_k, 2m_k)&\ge&2N-\frac{e(M_k)}{2},\lb{2.25}\\
i(\gamma_k, 2m_k)+\nu(\gamma_k, 2m_k)&\le&2N+\frac{e(M_k)}{2},\lb{2.26}
\eea
for every $k=1,\ldots,q$, where $S^+_{M_k}(1)$ is the splitting number of $M_k$.

More precisely, by (4.10) and (4.40) in \cite{LoZ} , we have
\bea m_k=\left(\left[\frac{N}{M\hat i(\gamma_k, 1)}\right]+\chi_k\right)M,\quad 1\le k\le q,\lb{2.27}\eea
where $\chi_k=0$ or $1$ for $1\le k\le q$ and $\frac{m_k\theta}{\pi}\in\Z$
whenever $e^{\sqrt{-1}\theta}\in\sigma(M_k)$ and $\frac{\theta}{\pi}\in\Q$
for some $1\le k\le q$.  Furthermore, given $M_0\in\N$,
by the proof of Theorem 4.1 of \cite{LoZ}, we may
further require  $M_0|N$ (since the closure of the
set $\{\{Nv\}: N\in\N, \;M_0|N\}$ is still a closed
additive subgroup of $\bf T^h$ for some $h\in\N$,
where we use notations  as (4.21) in \cite{LoZ}.
Then we can use the proof of Step 2 in Theorem 4.1
of \cite{LoZ} to get $N$).}

\setcounter{figure}{0}
\setcounter{equation}{0}
\section{Proof of Theorem 1.1}%Section 3

Firstly we make the following assumption

{\bf (FCG)} {\it Suppose that there exist only finitely many closed geodesics $c_k$ $(k=1, \cdots, q)$ on $(S^n,F)$ with reversibility $\lambda$ and flag curvature $K$ satisfying $\left(\frac{\lm}{1+\lm}\right)^2<K\le 1$.}

\medskip

If the flag curvature $K$ of $(S^n, F)$ satisfies
$\left(\frac{\lambda}{\lambda+1}\right)^2<K\le 1$,
then every non-constant closed geodesic $c$ must satisfy
\bea i(c)&\ge& n-1, \lb{3.1}\\
\hat i(c)&>& n-1, \lb{3.2}\eea
where (\ref{3.1}) follows from Theorem 3 and Lemma 3 of \cite{Rad3},
(\ref{3.2}) follows from Lemma 2 of \cite{Rad4}.
Thus it follows from Theorem 2.2 of \cite{LoZ} (or, Theorem 10.2.3 of \cite{Lon3}) that
\bea i(c^{m+1})-i(c^m)-\nu(c^m)\ge i(c)-\frac{e(P_c)}{2}\ge 0,\quad\forall m\in\N.\lb{3.3}\eea
Here the last inequality holds by (\ref{3.1}) and the fact that $e(P_c)\le 2(n-1)$.

It follows from (\ref{3.2}) and Theorem 2.7 that there exist infinitely many $(q+1)$-tuples $(N, m_1,\cdots, m_q)\in\N^{q+1}$ such that for any $1\le k\le q$, there holds
\bea i(c_k^{2m_k -1})+\nu(c_k^{2m_k -1})&=&
2N-\left(i(c_k)+2S^+_{M_k}(1)-\nu(c_k)\right), \lb{3.4}\\
i(c_k^{2m_k})&\ge& 2N-\frac{e(P_{c_k})}{2},\lb{3.5}\\
i(c_k^{2m_k})+\nu(c_k^{2m_k})&\le& 2N+\frac{e(P_{c_k})}{2},\lb{3.6}\\
i(c_k^{2m_k+1})&=&2N+i(c_k).\lb{3.7}\eea
Note that by List 9.1.12 of \cite{Lon3} and the fact $\nu(c_k)=p_{k_-}+2p_{k_0}+p_{k_+}$, we obtain
\bea 2S^+_{M_k}(1)-\nu(c_k) =2(p_{k_-}+p_{k_0}) - (p_{k_-}+2p_{k_0}+p_{k_+})=p_{k_-} -p_{k_+}.\lb{3.8}\eea

So by (\ref{3.1})-(\ref{3.8}) and the fact $e(P_{c_k})\le 2(n-1)$ it yields
\bea i(c_k^{m})+\nu(c_k^m)&\le& 2N-i(c_k)-p_{k_-} +p_{k_+}, \qquad\forall\ 1\le m<2m_k,\lb{3.9}\\
i(c_k^{2m_k})+\nu(c_k^{2m_k})&\le& 2N+\frac{e(P_{c_k})}{2}\le 2N+n-1,\lb{3.10}\\
2N+n-1&\le &i(c_k^{m}),\qquad\forall\ m>2m_k.\lb{3.11}\eea

In addition, the precise formulae of $i(c_k^{2m_k})$ and  $i(c_k^{2m_k})+\nu(c_k^{2m_k})$ can be computed as follows (cf. (3.16) and (3.21) of \cite{Dua} for the details)
\bea i(c_k^{2m_k})
&=& 2N -S^+_{M_k}(1)-C(M_k)+2\Delta_k,\lb{3.12}\\
i(c_k^{2m_k})+\nu(c_k^{2m_k})
    &=&2N+p_{k_0}+p_{k_+}+q_{k_-}+q_{k_0}\nn\\
      &&+2r_{k_0}'-2(r_{k_\ast}-r_{k_\ast}')+2r_k'-r_k+2\Delta_k,\quad k=1,\cdots,q, \lb{3.13}\eea
where $r_k',r_{k_\ast}'$ and $r_{k_0}'$ denote the number of normal forms $R(\th), N_2(e^{\aa\sqrt{-1}},A)$ and $N_2(e^{\bb\sqrt{-1}},B)$ with $\th,\aa,\bb$ being the rational multiple of $\pi$ in (\ref{2.16}) of Theorem 2.6 respectively, and
\bea \Delta_k \equiv \sum_{0<\{m_k\th/\pi\}<\delta}S^-_{M_k}(e^{\sqrt{-1}\th})\le r_k-r_k'+r_{k_\ast}-r_{k_\ast}',\quad
C(M_k) \equiv \sum_{\th\in(0,2\pi)}S^-_{M_k}(e^{\sqrt{-1}\th}), \lb{3.14}\eea
where $\delta>0$ is a small enough number (cf. (4.43) of \cite{LoZ}) and the estimate of $\Delta_k$ follows from the inequality (3.18) of \cite{Dua}.

Under the assumption ({\bf FCG}), Theorem 1.1 of \cite{Dua} shows that there exist at least two elliptic closed geodesics $c_1$ and $c_2$ on $(S^n,F)$
whose flag curvature satisfies $\left(\frac{\lm}{1+\lm}\right)^2<K\le 1$. Next Lemma lists some properties of these two closed geodesics which will be
useful in the proof of Theorem 1.1.

\medskip

{\bf Lemma 3.1.} (cf. Section 3 of \cite{Dua}) {\it Under the assumption (FCG), there exist at least two elliptic closed geodesics $c_1$ and $c_2$ on $(S^n,F), n\ge 2$
whose flag curvature satisfies $\left(\frac{\lm}{1+\lm}\right)^2<K\le 1$. Moreover, there exist infinitely many pairs of $(q+1)$-tuples $(N, m_1, m_2,\cdots, m_q)\in\N^{q+1}$ and $(N', m_1', m_2',\cdots, m_q')\in\N^{q+1}$ such that
\bea && i(c_1^{2m_1})+\nu(c_1^{2m_1})=2N+n-1,\quad \ol{C}_{2N+n-1}(E,c_1^{2m_1})=\Q,\lb{3.15}\\
&& i(c_2^{2m_2'})+\nu(c_2^{2m_2'})=2N'+n-1,\quad \ol{C}_{2N'+n-1}(E,c_2^{2m_2'})=\Q,\lb{3.16}\\
&& p_{k_-}=q_{k_+}=r_{k_\ast}=r_{k_0}-r_{k_0}'=h_k=0,\quad k=1,2,\lb{3.17}\\
&& r_1-r_1'=\Delta_1\ge 1,\qquad r_2-r_2'=\Delta_2'\ge 1, \lb{3.18}\\
&& \Delta_k+\Delta_k'=r_k-r_k',\quad k=1,2,\lb{3.19}\eea
where we can require $(n-1)|N$ or $(n-1)|N'$ as remarked in Theorem 2.7 and \bea \Delta_k' \equiv \sum_{0<\{m_k'\th/\pi\}<\delta}S^-_{M_k}(e^{\sqrt{-1}\th}),\quad k=1,2.\lb{3.20}\eea}

{\bf Proof.} In fact, all these properties have already been obtained in Section 3 of \cite{Dua} and here we only list references.
More precisely, (\ref{3.15}) follows from Claim 1 and arguments between (3.25) and (3.26) in \cite{Dua}. (\ref{3.16}) follows from Claim 3 and
similar arguments between (3.25) and (3.26) as those of $c_1$ in \cite{Dua}. (\ref{3.17}) and (\ref{3.18}) follow from (3.25), Claim 2 and Claim 3 in \cite{Dua}.
(\ref{3.19}) follows from (3.31) of \cite{Dua} and (\ref{3.17}). In a word, the properties of $c_1$ and $c_2$ is symmetric. \hfill\hb

\medskip

{\bf Lemma 3.2.} {\it Under the assumption (FCG), for these two elliptic closed geodesics $c_1$ and $c_2$ found in Lemma 3.1, the following further properties hold:
\bea &&i(c_k^{m})+\nu(c_k^m)\le 2N-i(c_k)+p_{k_+}\le 2N-1, \quad\forall\ 1\le m<2m_k, k=1,2,\lb{3.21}\\
 &&i(c_2^{2m_2})+\nu(c_2^{2m_2})\le 2N+n-3,\lb{3.22}\\
&&i(c_1^{2m_1'})+\nu(c_1^{2m_1'})\le 2N'+n-3,\lb{3.23}\eea
where the equalities in (\ref{3.22}) and (\ref{3.23}) hold respectively if and only if
\bea p_{k_0}+p_{k+}+q_{k_-}+q_{k_0}+2r_{k_0}'+r_k'=n-2,\quad r_k-r_k'=1,\quad k=1,2.\lb{3.24}\eea}

{\bf Proof.} Here we only give the proof of $c_2$. And the proof of $c_1$ is the same by using some information of $N'$ and $2m_1'$ instead of those of $N$ and $2m_2$
in the following arguments.

In fact, by (\ref{3.18}) and (\ref{3.19}), there holds $\Delta_2=0$. Then, together with the fact $r_{2_\ast}=0$ from (\ref{3.17}), it follows from (\ref{3.13}) that
\bea i(c_2^{2m_2})+\nu(c_2^{2m_2})&=&2N+(p_{2_0}+p_{2+}+q_{2_-}+q_{{2}_{0}}+2r_{2_0}'+r_2')-(r_2-r_2').\lb{3.25}\eea

Note that by (\ref{2.17}) we have \bea p_{2_0}+p_{2+}+q_{2_-}+q_{{2}_{0}}+2r_{2_0}'+r_2'+(r_2-r_2')=p_{2_0}+p_{2+}+q_{2_-}+q_{{2}_{0}}+2r_{2_0}'+r_2\le n-1.\lb{3.26}\eea
Therefore by (\ref{3.18}) we get
\bea p_{2_0}+p_{2+}+q_{2_-}+q_{{2}_{0}}+2r_{2_0}'+r_2'\le n-1-(r_2-r_2')\le n-2,\lb{3.27}\eea
which, together with (\ref{3.25}), yields
\bea i(c_2^{2m_2})+\nu(c_2^{2m_2})&=&2N+(p_{2_0}+p_{2+}+q_{2_-}+q_{{2}_{0}}+2r_{2_0}'+r_2')-(r_2-r_2')\nn\\
  &\le& 2N+n-3,\lb{3.28}\eea
where the equality holds if and only if
\bea p_{2_0}+p_{2+}+q_{2_-}+q_{{2}_{0}}+2r_{2_0}'+r_2'=n-2,\quad r_2-r_2'=1.\lb{3.29}\eea
So (\ref{3.22})) and (\ref{3.23}) follow from (\ref{3.28}) and (\ref{3.29}). Then (\ref{3.21}) follows from (\ref{3.1}) and (\ref{3.27}).

This completes the proof of Lemma 3.2.  \hfill\hb

\medskip

{\bf Lemma 3.3.} {\it Under the assumption (FCG), for these two elliptic closed geodesics $c_1$ and $c_2$ found in Lemma 3.1, there holds
\bea k_{\nu(c_k^{n(c_k)})}^{\ep(c_k^{n(c_k)})}(c_k^{n(c_k)})=1,\quad k_{j}^{\ep(c_k^{n(c_k)})}(c_k^{n(c_k)})=0,\quad\forall\  0\le j<\nu(c_k^{n(c_k)}),\ k=1,2. \lb{3.30}\eea}

{\bf Proof.} We only give the proof for $c_1$. The proof for $c_2$ is the same.

Firstly, by (\ref{3.15}) and Lemma 2.1, we have
\bea 1&=&\dim \ol{C}_{2N+n-1}(E,c_1^{2m_1})\nn\\
     &=&\dim H_{{2N+n-1}-i(c_1^{2m_1})}(N_{c_1^{2m_1}}\cup\{c_1^{2m_1}\},N_{c_1^{2m_1}})^{\ep(c_1^{2m_1})\Z_{2m_1}}\nn\\
     &=&\dim H_{\nu(c_1^{2m_1})}(N_{c_1^{2m_1}}\cup\{c_1^{2m_1}\},N_{c_1^{2m_1}})^{\ep(c_1^{2m_1})\Z_{2m_1}}\nn\\
     &=&k_{\nu(c_1^{2m_1})}^{\ep(c_1^{2m_1})}(c_1^{2m_1}),\eea
which implies that $k_{j}^{\ep(c_1^{2m_1})}(c_1^{2m_1})=0$, for any $0\le j<\nu(c_1^{2m_1})$ by (ii) of Lemma 2.3.
In addition, note that $n(c_1)|2m_1$ and $\nu(c_1^{2m_1})=\nu(c_1^{n(c_1)})$ by (\ref{2.7}) and (\ref{2.27}), there holds $k_{j}^{\ep(c_1^{2m_1})}(c_1^{2m_1})=k_{j}^{\ep(c_1^{n(c_1)})}(c_1^{n(c_1)})$ for any $0\le j\le\nu(c_1^{2m_1})$ by (iv) of Lemma 2.2. Thus (\ref{3.30}) holds. \hfill\hb

\medskip

{\bf Proof of Theorem 1.1:}

\medskip

In order to prove Theorem 1.1, we make the following assumption

\smallskip

{\bf (TCG)} {\it Suppose that there exist exactly two elliptic closed geodesics $c_1$ and $c_2$ on $(S^n,F)$ with reversibility $\lambda$ and flag curvature $K$ satisfying $\left(\frac{\lm}{1+\lm}\right)^2<K\le 1$.}

\smallskip

{\bf Claim 1:} {\it Under the assumption (TCG), $M_{2N+n-3}$ only can be contributed by $c_2^{2m_2}$, i.e.,
\bea M_{2N+n-3} =\sum_{1\le j\le 2,\; m\ge 1}\dim{\ol{C}}_{2N+n-3}(E, c^m_j)=\dim{\ol{C}}_{2N+n-3}(E, c_2^{2m_2}),\lb{3.32}\eea and
\bea &&i(c_2^{2m_2})+\nu(c_2^{2m_2})=2N+n-3,\lb{3.33}\\
&&i(c_1^{2m_1'})+\nu(c_1^{2m_1'})=2N'+n-3,\lb{3.34}\\
&&p_{k_0}+p_{k+}+q_{k_-}+q_{k_0}+2r_{k_0}'+r_k'=n-2,\quad r_k-r_k'=1,\quad k=1,2.\lb{3.35}\eea}
\indent In fact, note that $i(c_k^{m})+\nu(c_k^{m})\le 2N-1$ by (\ref{3.21}), this implies that all iterates $c_k^{m}$ with $1\le m\le 2m_k-1$, $k=1,2$ have no contribution to $M_{q},\ q\ge 2N$ by (i) of Lemma 2.2. On the other hand, by (\ref{3.11}) we know that all iterates $c_k^{m}$ with $m\ge 2m_k+1$, $k=1,2$ have no contribution to $M_{q},\ q\le 2N+n-2$ by (i) of Lemma 2.2. By (\ref{3.30}), the iterate $c_1^{2m_1}$ only contributes $1$ to $M_{2N+n-1}$ and has no contribution to $M_q$ with $q\neq 2N+n-1$.

If $i(c_2^{2m_2})+\nu(c_2^{2m_2})<2N+n-3$ by (\ref{3.22}), then $c_2^{2m_2}$ has no contradiction to $M_{2N+n-3}$ by Lemma 2.2. So there holds $M_{2N+n-3}=0$. However, it yields $b_{2N+n-3}\ge 1$ by Lemma 2.4. This gives a contradiction $0=M_{2N+n-3}\ge b_{2N+n-3}\ge 1$ by Theorem 2.5. Thus $i(c_2^{2m_2})+\nu(c_2^{2m_2})=2N+n-3$ and (\ref{3.32}) holds.

Similarly, (\ref{3.34}) can be obtained by using $N'$ and $2m_1'$ instead of $N$ and $2m_2$. Then, together with (\ref{3.24}), this completes the proof of Claim 1.

\medskip

{\bf Claim 2:} {\it Under the assumption (TCG), there exists at least one closed geodesic, without loss of generality, saying $c_1$, satisfying $p_{1_+}=n-2$ and $i(c_1)=n-1$.
And then there must be $i(c_2)=n-1$ and the analytic period defined by (\ref{2.7}) satisfies $n(c_2)\neq 1$.}

\medskip

In fact, if $i(c_k)>n-1$ or $p_{k_+}<n-2$ by (\ref{3.1}) and (\ref{3.27}), then by (\ref{3.17}), (\ref{3.9}) becomes
\bea i(c_k^{m})+\nu(c_k^{m})\le 2N-i(c_k)+p_{k_+}\le 2N-2,\quad 1\le m<2m_k,\ k=1,2.\lb{3.36}\eea
\indent Thus it yields $M_{2N-1}=0$ and $M_{2N+n-3}=1$ by (\ref{3.11}) and Claim 1.

When $n=3$, we obtain a contradiction $1=M_{2N}\ge b_{2N}=2$ by Lemma 2.4 and Theorem 2.5.

When $n\ge 4$, we obtain a contradiction $0=M_{2N-1}+M_{2N}\ge b_{2N-1}+b_{2N}\ge 1$ by Lemma 2.4 and Theorem 2.5.

So, without loss of generality, there holds
$p_{1_+}=n-2$ and $i(c_1)=n-1$. Then the analytic period $n(c_1)=1$ by (\ref{2.7}), which, together with Lemma 3.3 and $i(c_1^m)\ge i(c_1)\ge n-1$, yields
\bea \ol{C}_{n-1}(E,c_1^m)=k_{n-1-i(c_1^m)}^{\ep(c_1^m)}(c_1^m) = 0,\quad m\ge 1,\nn\eea
which shows that $c_1^m$ has no contribution to non-zero $M_{n-1}(\ge b_{n-1}=1)$. Thus there must be $i(c_2)=n-1$ and $n(c_2)>1$. In fact, if $n(c_2)=1$ or $i(c_2)>n-1$, then it can be shown that $c_2^m$ has no contribution to $M_{n-1}$ by Lemma 3.3 or Lemma 2.2. This contradiction completes this proof.

\medskip

Next we will carry out the proof according to the value of $n$ for $S^n$ in four cases.

\medskip

{\bf Case 1:} {\it $n=3$.}

\medskip

The existence of at least three closed geodesics in this case has been proved in \cite{Wan1}. Here we give a new and more simple proof.

In this case, there holds and $b_{2N}=2$ by Lemma 2.4.
In addition, it follows from (\ref{3.11}), (\ref{3.21}) and Claim 1 that among all iterates $c_k^m$, $m\ge 1$ of $c_k$, $k=1,2$, only $c_2^{2m_2}$ contributes $1$ to $M_{2N}$. So, under the assumption ({\bf TCG}), there holds $M_{2N}=1$, which contradicts to the Morse inequality $1=M_{2N}\ge b_{2N}=2$ and completes the proof in the case of $n=3$.

\medskip

{\bf Case 2:} {\it $n=4$.}

\medskip

According to (\ref{3.35}) and Claim 2, there
exists a continuous path $f\in\Om^0(P_{c_1})$ such that $f(0)=P_{c_1}$ and $f(1)=N_1(1,-1)^{\dm 2}\dm R(\th_{1}), \frac{\th_1}{\pi}\notin \Q$. Then by Theorem 2.6 we obtain
\bea i(c_1)=3,\quad i(c_1^m)&=&2m+2 E\left(\frac{m\th_1}{2\pi}\right) - 1,\quad \nu(c^m)=2, \quad\forall m\ge1. \lb{3.37}\eea

So, for the closed geodesic $c_1$, by Lemma 3.3 and (\ref{3.37}) we have
\bea \hat{\chi}(c_1)=\sum_{l=0}^{2}(-1)^{i(c_1)+l}k_l^+(c_1)=-k_0(c_1)+k_1^+(c_1)-k_2^+(c_1)=-1.\lb{3.38}\eea

Note that $2\times 2$ identity matrix $I_2$ and $-I_2$ can be viewed as a rotation matrix $R(\th)$ with $\th=2\pi$ and $\th=\pi$, respectively. Since $p_{2_0}+p_{2+}+q_{2_-}+q_{2_0}+2r_{2_0}'+r_2'=2$ by (\ref{3.35}), we only consider $p_{2+}+q_{2_-}+2r_{2_0}'+r_2'=2$ and $r_2-r_2'=1$. By Claim 2, $n(c_2)\neq 1$ implies $p_{2+}\neq 2$. Therefore \bea p_{2+}+q_{2_-}+2r_{2_0}'+r_2'=2,\quad p_{2+}\neq 2,\lb{3.39}\eea
which yields
\bea 4\ge\nu(c_2^{2m_2})=\nu(c_2^{n(c_2)})=p_{2+}+q_{2_-}+2r_{2_0}'+2r_2'\ge 2.\lb{3.40}\eea

If $\nu(c_2^{n(c_2)})=4$, there must be $r_2'=2$ by (\ref{3.40}). In this case, there holds $i(c_2^m)\in 2\N-1$, $\forall m\ge 1$ by (4.7) and (4.8) of \cite{BaL} and the symplectic additivity (cf. Theorem 9.1.10 of \cite{Lon3}), and $\nu(c_2^m)\le 2$, $\forall 1\le m<n(c_2)$ since $n(c_2)>1$ by Claim 2.

If $\nu(c_2^{n(c_2)})=3$, there holds either $r_2'=1$ and $p_{2_+}=1$, or $r_2'=1$ and $q_{2_-}=1$ by (\ref{3.40}). In either case, by (4.3), (4.6)-(4.8) of \cite{BaL} and the symplectic additivity,  for any $1\le m<n(c_2)$, either there holds $\nu(c_2^m)\le 1$, or there holds $\nu(c_2^m)=2$ and $i(c_2^m)\in 2\N-1$ which only happens in the case of $r_2'=1$ and $q_{2_-}=1$.

If $\nu(c_2^{n(c_2)})=2$, there must be $r_{2_0}'=1$ or $q_{2_-}=2$ by (\ref{3.40}). In this case, $n(c_2)>1$ and $\nu(c_2^m)=0$, $\forall\ 1\le m<n(c_2)$.

In summary, only one of the following two cases can happen:

\smallskip

(i) $\nu(c_2^m)\le 1$ for some $1\le m<n(c_2)$,

(ii) $\nu(c_2^m)=2, i(c_2^m)+\nu(c_2^m)\in 2\N-1$ for some $1\le m<n(c_2)$.

\smallskip

Firstly, by Lemma 3.3 and (\ref{3.33}) we have
\bea  {\chi}(c_2^{n(c_2)})&=&\sum_{l=0}^{\nu(c_2^{n(c_2)})} (-1)^{i(c_2^{n(c_2)})+l}k_l^{\ep(c_2^{n(c_2)})}(c_2^{n(c_2)})=-k_{\nu(c_2^{n(c_2)})}^{\ep(c_2^{n(c_2)})}(c_2^{n(c_2)})=-1,\lb{3.41}\eea

If (i) happens, then for some $1\le m<n(c_2)$ satisfying $\nu(c_2^m)\le 1$, by (ii) of Lemma 2.2, it yields
\bea {\chi}(c_2^{m})&=&\sum_{l=0}^{\nu(c_2^m)} (-1)^{i(c_2^{m})+l}k_l^{\ep(c_2^m)}(c_2^{m})=(-1)^{i(c_2^{m})}\left(k_0^{\ep(c_2^m)}(c_2^m)-k_1^{\ep(c_2^m)}(c_2^m)\right)
  \ge -1.\lb{3.42}\eea

If (ii) happens, then for some $1\le m<n(c_2)$ satisfying $\nu(c_2^m)=2$, by (ii) of Lemma 2.2, it yields
\bea {\chi}(c_2^{m})&=&\sum_{l=0}^2 (-1)^{i(c_2^{m})+l}k_l^{\ep(c_2^m)}(c_2^{m})=-k_0(c_2^m)+k_1^+(c_2^m)-k_2^+(c_2^m)\ge -1.\lb{3.43}\eea

Now by (\ref{3.41})-(\ref{3.43}), we obtain
\bea \hat{\chi}(c_2)=\frac{1}{n(c_2)}\sum_{m=1}^{n(c_2)}\chi(c_2^m)\ge -1.\lb{3.44}\eea

Note that $\hat{i}(c_k)>3$, $k=1,2$ by (\ref{3.2}), so it follows from (\ref{3.38}) and (\ref{3.44}) that
\bea \frac{\hat{\chi}(c_1)}{\hat{i}(c_1)}=-\frac{1}{\hat{i}(c_1)}>-\frac{1}{3},\qquad
\frac{\hat{\chi}(c_2)}{\hat{i}(c_2)}\ge-\frac{1}{\hat{i}(c_2)}>-\frac{1}{3},\eea
which, together with Theorem 2.3, yields
\bea -\frac{2}{3}=\frac{\hat{\chi}(c_1)}{\hat{i}(c_1)}+\frac{\hat{\chi}(c_2)}{\hat{i}(c_2)}>-\frac{2}{3}.\eea
This contradiction completes the proof of Theorem 1.1 in the Case of $n=4$.

\medskip

{\bf Case 3:} {\it $n=5$.}

\medskip

In this case, note that $(n-1)|N$, there holds $b_{2N}=2$ by Lemma 2.4.
In addition, it follows from (\ref{3.11}), (\ref{3.21}) and Claim 1 that all iterates $c_k^m$, $m\ge 1$ of $c_k$, $k=1,2$ have no contribution to $M_{2N}$. So, under the assumption ({\bf TCG}), there holds $M_{2N}=0$, which contradicts to the Morse inequality $0=M_{2N}\ge b_{2N}=2$ and completes the proof in the case of $n=5$.

\medskip

{\bf Case 4:} {\it $n\ge 6$.}

\medskip

In this case, note that $(n-1)|N$, there holds $b_{2N+n-5}=1$ by Lemma 2.4. In addition, it follows from (\ref{3.11}), (\ref{3.21}) and Claim 1 that all iterates $c_k^m$, $m\ge 1$ of $c_k$, $k=1,2$ have no contribution to $M_{2N+n-5}$.

So, under the assumption ({\bf TCG}), there holds $M_{2N+n-5}=0$, which contradicts to the Morse inequality $0=M_{2N+n-5}\ge b_{2N+n-5}=1$. Thus there must exist the third closed geodesic $c_3$ such that the iterates $c_3^m$, $m\ge 1$ have contributed at least $1$ to $M_{2N+n-5}$. On the other hand, it follows from (\ref{3.1}), (\ref{3.9}) and (\ref{3.11}) that $M_{2N+n-5}$ can be contributed only by the iterate $c_3^{2m_3}$. So we have
\bea \ol{C}_{2N+n-5}(E,c_3^{2m_3})\neq 0. \lb{3.47}\eea

If $c_3$ is a hyperbolic closed geodesic which implies that $S^+_{M_3}(1)=C(M_3)=\Delta_3=0$, then by (\ref{3.13}) it yields $i(c_3^{2m_3})+\nu(c_3^{2m_3})=2N$. Then by Lemma 2.2 we obtain $\ol{C}_{2N+n-5}(E,c_3^{2m_3})=0$ since $2N+n-5\ge 2N+1$ in this case. This contradiction with (\ref{3.47}) shows that the closed geodesic $c_3$ must be non-hyperbolic.  \hfill\hb

\medskip

{\bf Acknowledgment.} The author sincerely thank the referee for her/his careful reading, valuable comments
and suggestions on this paper.

\bibliographystyle{abbrv}

\end{document}